\RequirePackage{fix-cm}
\documentclass{article}

\usepackage{arxiv}

\usepackage[utf8]{inputenc} % allow utf-8 input
\usepackage[T1]{fontenc}    % use 8-bit T1 fonts
\usepackage[colorlinks=true,urlcolor = black,citecolor=black,linkcolor=blue]{hyperref}
\usepackage{url}            % simple URL typesetting
\usepackage{booktabs}       % professional-quality tables
\usepackage{nicefrac}       % compact symbols for 1/2, etc.

\usepackage{lipsum}         % Can be removed after putting your text content
\usepackage{graphicx}
\usepackage{doi}

\usepackage{color}
\usepackage{amssymb,amsfonts,amsmath,amsthm,amsbsy} % Math symbols.
\usepackage{tabularx} % Helpful for full text width tables.
\usepackage{mathtools,mathrsfs}
\usepackage{mhchem}
\usepackage{enumerate}
\usepackage{bm}
\usepackage{bbm}
\usepackage{multirow}
\usepackage{threeparttable}
\usepackage{float}

\usepackage[table,dvipsnames, svgnames]{xcolor}

\usepackage[bottom]{footmisc}
\usepackage[numbers,sort&compress]{natbib}

\usepackage{indentfirst}
\usepackage{etoolbox,siunitx}

\usepackage{afterpage}
\robustify\bfseries
\usepackage{adjustbox}
\usepackage{chngpage}
\usepackage{rotating}
\usepackage{hvfloat}
\usepackage{soul}

\usepackage{subcaption}
\usepackage{caption}

\usepackage[italicComments=true,indLines=true,noEnd=false]{algpseudocodex}
\usepackage[english]{babel}
\hypersetup{colorlinks=true, citecolor=blue}
\usepackage{stackrel}
\usepackage[nopatch=eqnum]{microtype}      % microtypography

\definecolor{abstract_background}{RGB}{235,235,235}%
\usepackage[title]{appendix}
\usepackage[capitalise,nameinlink]{cleveref}

\title{Timescale dynamics of COVID-19 pandemic waves: The case of Greece}

\author{
\textbf{\href{https://orcid.org/0000-0002-4788-9877}{\includegraphics[scale=0.06]{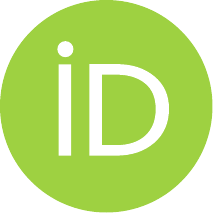}\hspace{1mm}Dimitris M. Manias}\textcolor{teal}{$^{1,}$}\thanks{Corresponding Author: dimitris.manias@ku.ac.ae}~, 
\href{https://orcid.org/0000-0001-9840-0018}{\includegraphics[scale=0.06]{orcid.pdf}\hspace{1mm}Dimitris G. Patsatzis}\textcolor{teal}{$^{2}$}, 
\href{https://orcid.org/0000-0002-9674-9491}{\includegraphics[scale=0.06]{orcid.pdf}\hspace{1mm}Dimitris A. Goussis}\textcolor{teal}{$^{3,4}$}}
{}\\[5pt]
\textcolor{teal}{$^{(1)}$} Department of Mathematics, \emph{Khalifa University of Science and Technology}, \\ Abu Dhabi 127788, United Arab Emirates \\
\textcolor{teal}{$^{(2)}$} Modelling Engineering Risk and Complexity, \emph{Scuola Superiore Meridionale}, Naples 80138, Italy \\
\textcolor{teal}{$^{(3)}$} Department of Mechanical Engineering, \emph{Khalifa University of Science and Technology}, \\ Abu Dhabi 127788, United Arab Emirates\\
\textcolor{teal}{$^{(4)}$} Research and Innovation Center on CO$_2$ and H$_2$,  \emph{Khalifa University of Science and Technology}, \\ Abu Dhabi 127788, United Arab Emirates
}

\renewcommand{\headeright}{}
\renewcommand{\undertitle}{}
\renewcommand{\shorttitle}{\textit{Time scale dynamics of COVID-19 pandemic waves: The case of Greece} }

\hypersetup{
pdftitle={Time scale dynamics of COVID-19 pandemic waves: The case of Greece },
pdfsubject={math.DS},
pdfauthor={Dimitris M. Manias},
pdfkeywords={Population dynamics, COVID-19, Computational Singular Perturbation, time scale analysis},
}

\begin{document}
\maketitle

\begin{abstract}
\colorbox{abstract_background}{\begin{minipage}{1\linewidth}
The results of an alternative methodology for making predictions about the COVID-19 pandemic in Greece are presented.~Instead of focusing on the various population profiles (subjected to instabilities introduced by the fitting process), this methodology focuses on the time scale that characterises the intensity and duration of the \emph{outbreak phase}.~Therefore, instead of predicting the peak of active cases, here their inflection point is predicted (the point where the increase of active cases stops accelerating and starts decelerating).~Since the inflection point precedes the peak, this methodology can serve as an early warning of the peak.~In addition, the paths between the various populations (healthy, exposed, infected, etc) that contribute the most to the \emph{outbreak phase} are identified.
\end{minipage}
}
\end{abstract}
\par\noindent\rule{\textwidth}{0.5pt}

% keywords can be removed

%%%%%%%%%%%%%%%%%%%%%%%%%%%%%%%%%%%%%%%
\section{Introduction}
\label{Intro}
The major effort in making predictions on the spread of COVID-19 focuses on the temporal evolution of various populations, notably active cases and deaths. Here, the employment of Computational Singular Perturbation (CSP) \citep{lam1994,lam1994csp,goussis1992} is proposed as an alternative methodological approach for the analysis of COVID-19 pandemics compartmental models.~The analysis is based on the time scale that characterises the \emph{outbreak phase} (the period in which the number of active cases increases) and employs tools used for the multi-scale analysis of models in chemical kinetics, biology and pharmacokinetics.%~The present work is based on experience acquired and tools developed during a KU-NTUA collaboration project, funded by KU, during May-December 2020.
CSP is an algorithmic method of Geometrical Singular Perturbation Analysis (GSPA), which is employed in multi-scale dynamical systems; i.e., systems that are driven by processes evolving in a wide range of fast/slow time scales.~In such systems, the processes characterized by the fast time scales become quickly exhausted and form low dimensional surfaces in the tangent space, well-known as Slow Invariant Manifolds (SIM) \citep{valorani2003csp,goussis2013role,maris2015hidden}.~The trajectories of the system are then confined to evolve on the SIM, governed by the processes which are characterized by the remaining slow time scales.~Recently, a physics-informed machine learning approach was introduced for the computation of SIMs, in the spirit of GSPA \cite{patsatzis2023slow}.~CSP offers a complete set of diagnostic tools for identifying the variables and physical processes related to (i) the fast time scales, which contribute to the formation of the SIM and (ii) the slow time scales that govern the dominant slow dynamics of the system.~The CSP tools have been employed for gaining physical understanding of chemical kinetic mechanisms \citep{lam1989,hadjinicolaou1998asymptotic, manias2016mechanism, rabbani2023comparative, rabbani2022dominant, rabbani2021ch3oh, manias2021nh3, manias2022mechanism, khalil2021no, khalil2019algorithmic, tingas2018ch4}
 and other multi-scale physical problems encountered in systems biology \citep{patsatzis2019new, patsatzis2016asymptotic} and population dynamics \citep{patsatzis2022algorithmic}.  

The utilization of CSP offers significant advantages for analyzing epidemiological models, particularly in compartmental contexts. CSP's relevance stems from its compatibility with multi-scale epidemiological models, as demonstrated in various studies \citep{gaucel2009using,brauer2019singular,souza2014multiscale,zhang2009singular,jardon2021geometric}. Unlike solution-based approaches, CSP focuses on system dynamics, avoiding issues like non-identifiability. It also provides a systems-level understanding, surpassing the capabilities of sensitivity analysis in distinguishing fast and slow dynamics. CSP's effectiveness in exploring dynamical systems with explosive eigenvalues, as seen in the initial exponential growth phase of COVID-19 \citep{chowell2011characterizing,wallinga2007generation,park2019practical}. The presence of a positive or explosive eigenvalue, which is a key characteristic during the exponential growth phase of COVID-19 pandemics, is a well-established concept \citep{patsatzis2021relation, barbarossa2020modeling}. CSP methodology, known for its effectiveness in studying dynamical systems with such eigenvalues \citep{valorani2020computational1, valorani2020computational2, manias2022effect, patsatzis2022algorithmic}, is therefore a logical choice for analyzing the rapid expansion of COVID-19.

\begin{figure}[h!]
\centering
\includegraphics[scale=0.35]{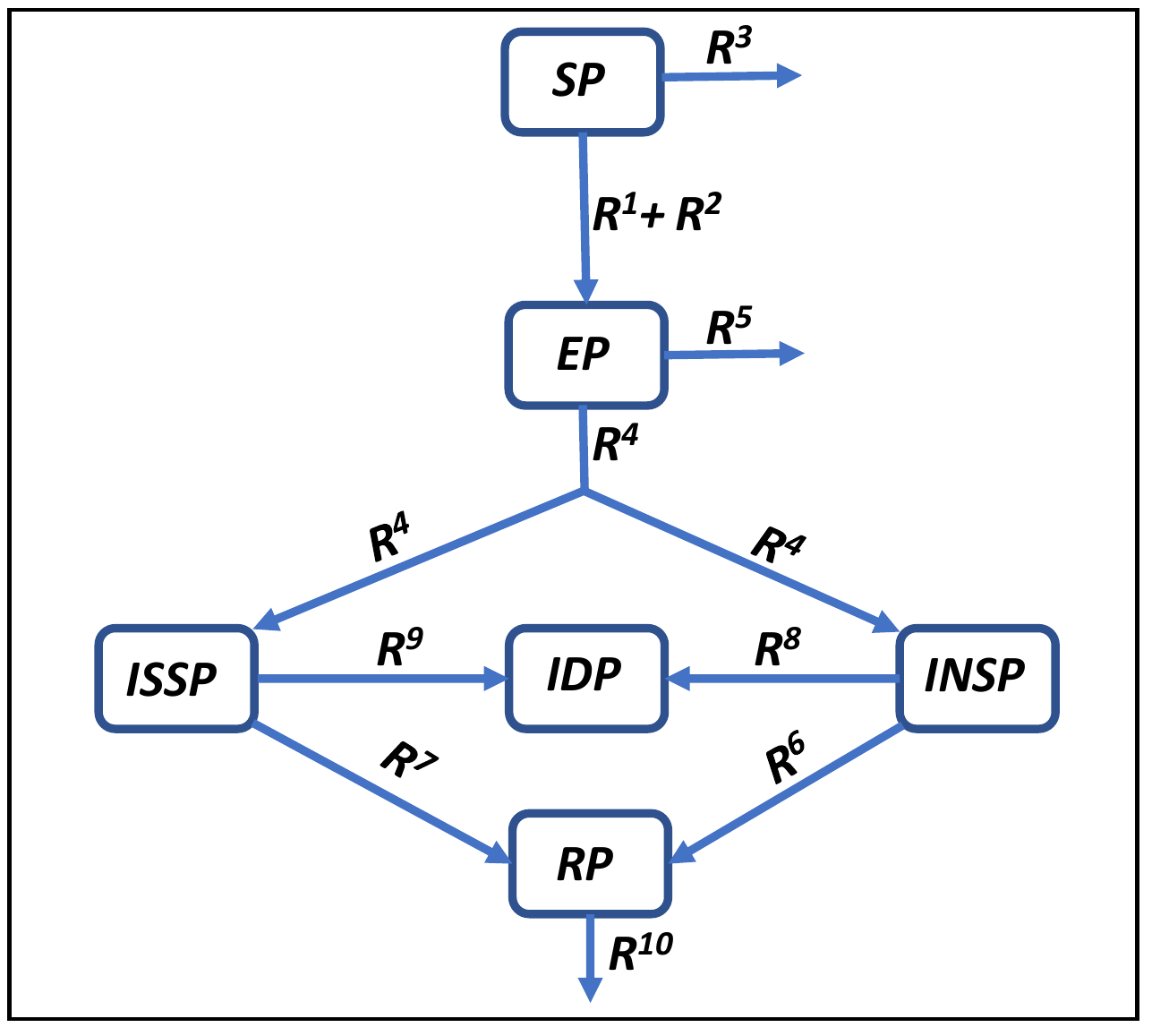}
\vspace{-3 pt}
\caption{The SEInsRD population model employed; SP: susceptible population (not infected), EP: exposed population (infected but not yet able to infect), INSP: infected normally susceptible population (able to infect), ISSP: infected severely susceptible population (able to infect), IDP: infected dead population, RP: recovered population.~$R^i$ denotes the transmission rate of the $i$-th path.}\label{fig:1x}
\end{figure}

In this study, we apply a novel approach to analyze COVID-19 pandemics in Greece using the SEInsRD model, shown in Fig. \ref{fig:1x}. We focus on the outbreak phases of the 4th, 5th, and 6th waves, identifying key populations and processes driving the pandemic dynamics using CSP diagnostic tools. Additionally, we explore the relationship between the system's time scales and the inflection point of the infected population. This approach allows early estimation of inflection points, aiding in predicting the plateau phase of the pandemic, as supported by previous studies \cite{fokas2020covid, fokas2020predictive}.

%---------------------------------------------------------------------------------------------------------------------------------------------
\section{Mathematical Model}
In this study, the SEInsRD population model was employed, a description of which is depicted in Fig. \ref{fig:1x}. Since research indicates that COVID-19 symptoms vary from mild to severe, the SEInsRD model is an adaptation of the SEIRD model, which expands the later by splitting the infected compartment into two sub-compartments: 'normally infected' (IN) and 'severely infected' (IS) individuals. This division reflects the immediate healthcare needs of severely infected individuals and leads to the modified SEInsRD model:

\begin{equation}
\dfrac{d}{dt} \begin{bmatrix} S \\ E \\ IN \\ IS \\ R \\ D \end{bmatrix}= \begin{bmatrix} - \beta_N S.IN - \beta_S S.IS - \mu_{TP} S \\  \beta_N S.IN + \beta_S S.IS -\sigma E - \mu_{TP} E \\ (1-ss) \sigma E - \gamma IN - \mu_N IN  \\ ss \sigma E - \gamma IS - \mu_S IS \\ \gamma (IN+IS) - \mu_{TP} R \\ \mu_N IN + \mu_S IS \end{bmatrix}
\end{equation}\label{eq:SEInsRD}
where the subscripts $N$ and $S$ indicate the normally and severely infected transmission $\beta$ and fatality $\mu$ ratios, $ss$ denotes the fraction of severely over normally infected individuals and $\mu_{TP}$ is the physiological death ratio.

The SEInsRD model can be rewritten to its - equal to density dependent - frequency dependent formulation, as:
\begin{align}
\dot{SP} & =-R^1-R^2-R^3 			&	\dot{EP} & =R^1+R^2-R^4-R^5 	\nonumber \\
\dot{INSP} & =(1-ss) R^4 - R^6 -R^8 	&	\dot{ISSP} & =ss R^4 - R^7 -R^9  \label{eq:SEI2RD}  \\
\dot{RP} & = R^6+R^7-R^{10}		 	&	\dot{IDP} & = R^8+R^9			\nonumber 
\end{align}

\noindent where $SP$, $EP$, $INSP$, $ISSP$, $RP$ and $IDP$ denotes the number of susceptible, exposed, infected normally symptomatic,  infected severely symptomatic,  recovered and  infected deceased individuals, respectively.~The rates $R^i$ are defined as:

\begin{align}
R^1 & =\dfrac{\beta_{INSP}}{TP}INSP~SP	&	R^2 & =\dfrac{\beta_{ISSP}}{TP}ISSP~SP &	R^3   & =\mu_{TP}SP 	\nonumber \\
R^4 & =aincp^{-1} EP					&	R^5 & =\mu_{TP}EP 					&	R^6   & =aip^{-1} INSP 	\nonumber \\
R^7 & =aip^{-1} ISSP 					&	R^8 & =\mu_{INSP}INSP 				&	R^9   & =\mu_{ISSP}ISSP  	\nonumber \\
R^{10} & =\mu_{TP}RP					&					 				&		\nonumber
\end{align}\label{eq:model}

\noindent where the transmission ratios are expressed in a frequency dependent representation as $\beta_N=\beta_{INSP}/TP$ and $\beta_S=\beta_{ISSP}/TP$, $\mu_N = \mu_{INSP}$, $\mu_S = \mu_{ISSP}$, $\gamma=aip^{-1}$ and $\sigma=aincp^{-1}$. In addition, Eq. (\ref{eq:SEI2RD}) yields the forming equilibration:
\begin{equation}
\dfrac{d(SP+EP+INSP+ISSP+RP+IDP)}{dt}=R^3+R^5+R^{10} =- \mu_{TP}(SP+EP+RP)
\end{equation}\label{eq:equil}

More details on the inflection points through analytic calculations are included in Appendix \ref{sec:appendix}.

%%%%%%%%%%%%%%%%%%%%%%%%%%%%%
%%%%
\section{Computational Singular Perturbation}
\label{sec:CSP}

The SEInsRD compartmental model of Eq. \eqref{eq:model} can be formed as an Ordinary Differential Equation (ODE) system:
\begin{equation}\label{eq:gov8}
\dfrac{d\mathbf{y}}{dt}= \mathbf{g(y)} =\sum_{k=1}^{K}\mathbf{S}_kR^k(\mathbf{y}) 
\end{equation}
where $\mathbf{y}$ is the $N$-dim. column \textit{state vector}, the elements of which are the number of individuals in the population group or their fraction over the population and $\mathbf{g(y)}$ is the $N$-dim. column \textit{vector field} that incorporates the  transition laws from a compartmental group to another.~In particular, the vector field $\mathbf{g(y)}$ is described by the $N$-dim. stoichiometric vector $\bm{S}_k$, which indicates the direction of the $K$ transitions from a group to another and $R_k$, which is the related transition rate/interaction; e.g., transmission, recovery, death rate.~According to CSP, the system in Eq. \eqref{eq:gov8} can be cast in the following form \cite{lam1989,lam1991}:
\begin{equation}
\label{eq:gov9}
\dfrac{d\mathbf{y}}{dt} = \mathbf{g}(\mathbf{y}) =  \sum_{n=1}^{N}\mathbf{a}_n(\mathbf{y})f^n(\mathbf{y}), \qquad f^n(\mathbf{y}) = \mathbf{b}^n(\mathbf{y}) \cdot \mathbf{g(y)} = \sum_{k=1}^{K} \left( \mathbf{b}^n(\mathbf{y}) \cdot \mathbf{S}_k\right)R^k(\mathbf{y}),
\end{equation} 
where $\mathbf{a}_n(\mathbf{y})$ and $\mathbf{b}^n(\mathbf{y})$ are the $N$-dim. CSP column and row, respectively, basis vectors of the $n$-th mode, which satisfy the orthogonality conditions $\mathbf{b}^i(\mathbf{y}) \cdot \mathbf{a}_j(\mathbf{y}) = \delta_j^i$ \cite{lam1989,lam1994}.~In this way, the vector field $\mathbf{g(y)}$ is decomposed in $n$ CSP modes, each having a different impact in the vector field and acting in a different timescale.~The impact of the $n$-th CSP mode is measured by the amplitude $f^n(\mathbf{y})$ (proper adjustment of the signs of the CSP basis vectors, set $f^n(\mathbf{y})$ always positive), which provides a measure of the projection of the vector field $\mathbf{g(y)}$ on the CSP vector $\mathbf{a}_n$.~When the system in Eq. \eqref{eq:gov9} exhibits $M$ timescales that are [i] of dissipative nature (i.e., when the components of the system that generate them tend to drive the system towards equilibrium) and [ii] much faster than the rest, the model can be reduced in:
\begin{equation}
\label{eq:gov10}
f^r(\mathbf{y})\approx 0\quad(r=1,\ldots,M), \qquad \dfrac{d\mathbf{y}}{dt} \approx \sum_{s=M+1}^{N}\mathbf{a}_s(\mathbf{y})f^s(\mathbf{y}),
\end{equation} 
the first relation of which is an $M$-dim. system of algebraic equations defining the manifold $\mathscr{M}$ (a low dimensional surface in phase-space, where the system is confined to evolve), while the second relation is an $N$-dim. system of ODEs governing the slow evolution of the system on this manifold.~The objective of asymptotic analysis techniques is to deliver the expressions in Eq.~\eqref{eq:gov10}, which are provided, order by order, by CSP in algorithmic fashion \cite{zagaris2004fast,kaper2015geometry}.~This is accomplished through the CSP vectors, $\mathbf{a}_i(\mathbf{y})$ and $\mathbf{b}^i(\mathbf{y})$ ($i=1,\ldots,N$), which can be approximated in leading order accuracy by the right and left, respectively, eigenvectors of the $N \times N$-dim.  Jacobian $\mathbf{J(y)}$ of $\mathbf{g(y)}$; i.e., $\mathbf{a}_i(\mathbf{y})=\boldsymbol{\alpha}_i(\mathbf{y})$ and $\mathbf{b}^i(\mathbf{y})=\boldsymbol{\beta}^i(\mathbf{y})$ \cite{lam1991,lam1989,tingas2016}.~In the following, for simplicity reasons, the dependency of $R^k$, $\mathbf{g}$, $\boldsymbol{\alpha}$, etc, to $\mathbf{y}$ will be removed, and be mentioned only when explicitly needed.

The interactions among the population groups have a different impact to the evolution of the system, either to the constraints along which the system is confined to evolve or to its slow evolution along them.~In particular, the $M$ constraints in Eq.~\eqref{eq:gov10} $f^r = (\boldsymbol{\beta}^r \cdot \mathbf{S}_1)R^1 + \ldots +  (\boldsymbol{\beta}^r \cdot \mathbf{S}_K)R^K \approx 0$ ($r=1,\ldots,M$) result from significant cancellations among some of the additive terms $(\boldsymbol{\beta}^r \cdot \mathbf{S}_k)R^k$ ($k=1,\ldots,K$).~In addition, the slow evolution of the system is characterized by the slow CSP modes, the impact of which is measured by the slow amplitudes $f^s$ in Eq.~\eqref{eq:gov10}.~The interactions that contribute significantly to (i) the formation of each of the $M$ constraints and (ii) the governing components of the slow system, are identified by the \textit{Amplitude Participation Index} (\textit{API}):
\begin{equation}
P^n_k = \frac{\displaystyle (\boldsymbol{\beta}^{n} \mathbf{S}_k)R^k}{\displaystyle  \sum\nolimits_{i=1}^{K}|    (\boldsymbol{\beta}^n \mathbf{S}_i)R^i | } \qquad (k=1,\ldots,K),
\label{eq:gov11}
\end{equation}
where by definition $\sum_{k=1}^{K}|P_k^n|=1$ \cite{lam1994,goussis2006,valorani2003}.~When employed to the fast $r=1,\ldots,M$ CSP modes, $P^r_k$ measures the relative contribution of the $k$-th interaction to the cancellations among the additive terms in $f^r\approx 0$, while when employed to the slow $s=M+1,\ldots,N$ CSP modes, $P^s_k$ measures the relative contribution of the $k$-th interaction to governing slow system.~Both $P^r_k$ and $P^s_k$ can be either positive or negative; the sum of positive and negative $P^r_k$ terms equaling 0.5, by definition \cite{manias2019topological, manias2019investigation, manias2019dynamics, manias2018analysis}.

The formation of the $M$ constraints is characterized by the $M$ fastest timescales, which by definition are of dissipative nature; see requirement [i] immediately before Eq.~\eqref{eq:gov10}.~On the contrary, the dynamics of the slow system in Eq.~\eqref{eq:gov10} is characterized either by the fastest of the $N-M$ dissipative slow timescales or by explosive timescales, when the later are present, since they tend to drive the system away from equilibrium.~The timescales are approximated by the inverse of the eigenvalues of the Jacobian $\mathbf{J}$, $\tau_n = |\lambda_n|^{-1}$ ($n=1,\ldots,N$).~The CSP diagnostic tool, {\em Timescale Participation Index} (\textit{TPI}), identifies the reactions significantly contributing to the generation of the timescales:
\begin{equation}\label{eq:gov12}
J^n_k = \dfrac{c_k^n}{\sum\nolimits_{i=1}^{K}|c^n_i|} \qquad (k=1,\ldots,K),
\end{equation}
where ${\lambda}_n = c_1^n + \ldots + c_{K}^n$ and by definition $\sum_{k=1}^{K}|J_k^n|=1$ \cite{goussis2005,goussis2006,diamantis2015}.~$c_k^n$ denotes the contribution of the $k$-th interaction to the $n$-th eigenvalue and can be calculated as $c_k^n=\boldsymbol{\beta}^n \nabla  \left(  \mathbf{S}_kR^k \right) \boldsymbol{\alpha}_n$, where $\sum_{k=1}^K \nabla \left(\mathbf{S}_kR^k\right)$ is the Jacobian $\mathbf{J}$.~$c_k^n$ can be either positive or negative and therefore, a negative (positive) $J^n_k$ implies that the $k$-th interaction contributes to a dissipative (explosive) character of the $n$-th timescale $\tau_n$ \cite{manias2015algorithmic, manias2015initiation}.~By definition, dissipative (explosive) timescales relate to the components of the system that tend to drive it towards (away from) equilibrium \cite{lam1989,lam1994}. TPI has been succesfuly used 

Each population group/compartment is associated differently to each CSP mode; e.g., a fast CSP mode is expected to be much more related to the infected compartment than the recovered one, since the former are playing a major role as opposed to the latter at the beginning of an epidemics.~The relation of the $m$-th CSP mode  ($m=1, \dots, M$) to the various population groups is identified by the \textit{Pointer} (\textit{Po}):
\begin{equation}\label{eq:gov13}
\mathbf{D}^m = diag \left[\boldsymbol{\alpha}_m\boldsymbol{\beta}^m \right] = \left[ \alpha^1_m\beta^m_1,\alpha^2_m\beta^m_2,\ldots,\alpha^{N}_m\beta^m_{N} \right] \qquad (m=1,\ldots,M),
\end{equation}
where, due to the orthogonality condition $\boldsymbol{\beta}^i \cdot \boldsymbol{\alpha}_j=\delta^i_j$, the sum of all $N$ elements of $\mathbf{D}^m$ equals unity, i.e. $\sum_{i=1}^N\alpha^i_m\beta^m_i=1$ \cite{goussis1992,lam1994,valorani2003,goussis2012}.~A relatively large value of $\alpha^i_m\beta^m_i$ indicates that the $i$-th population group is strongly associated to $m$-th CSP mode and the $m$-th timescale.~A value of $D^m_i$ close to unity suggests that the $i$-th population group is in \textit{Quasi Steady-State} (QSS) \cite{goussis2012}.

%-------------------------------------------------------------------------------------------------------------------------------------------------
\section{Results}

The fastest of the \emph{explosive time scales} in the dynamics of the process characterises the \emph{outbreak phase} from the start of the process, up to the point where these time scales cease to exist (this is the period in which the system accelerates at the start of the outbreak phase)  and that the remaining period of this phase is characterised by the fastest of the \emph{dissipative time scales} (this is the period in which the system decelerates as it approaches the end of the outbreak phase, known as the \emph{peak}).~The transition from the accelerating to the decelerating period, defines the \emph{inflection point} in the profile of the active cases, during the \emph{outbreak phase}.~The development of the \emph{inflection point}  is an early warning that the \emph{peak} will be reached soon.%~This feature was validated during the KU-NTUA project, by analysing outbreak phases in a large number of countries (USA, Canada, Spain, France, Italy, Austria, Germany, China, S. Korea, N. Zealand, etc).

The explosive and dissipative time scales that characterize the outbreak phase were computed by fitting the available data (from Worldometer database, number of exposed, infected and dead) to the SEInsRD model, shown in Fig.~\ref{fig:1x}.%~Work done during the KU-NTUA project, demonstrated that other population models (SIR, SEIR etc) provide qualitatively similar results.

%-------------------------------------------------------------------------------------------------------------------------------------------------
\subsection{The 4th wave}
\label{sec:4th_wave}
In this section, results for the 4th wave will be presented.~For validation purposes, the analysis was done by considering solutions obtained by fitting data from two periods in 2021:
\begin{enumerate}[i)]
   \item Period A, July 1 - July 8, 8 days long,
   \item Period B, July 1 - July 18, 18 days long,
\end{enumerate}
where July 1 marks the start of the 4th wave.

%\subsubsection{Results}
%\label{sec:4th_wave_results}

\begin{figure}[h!]
\centering
\includegraphics[scale=0.30]{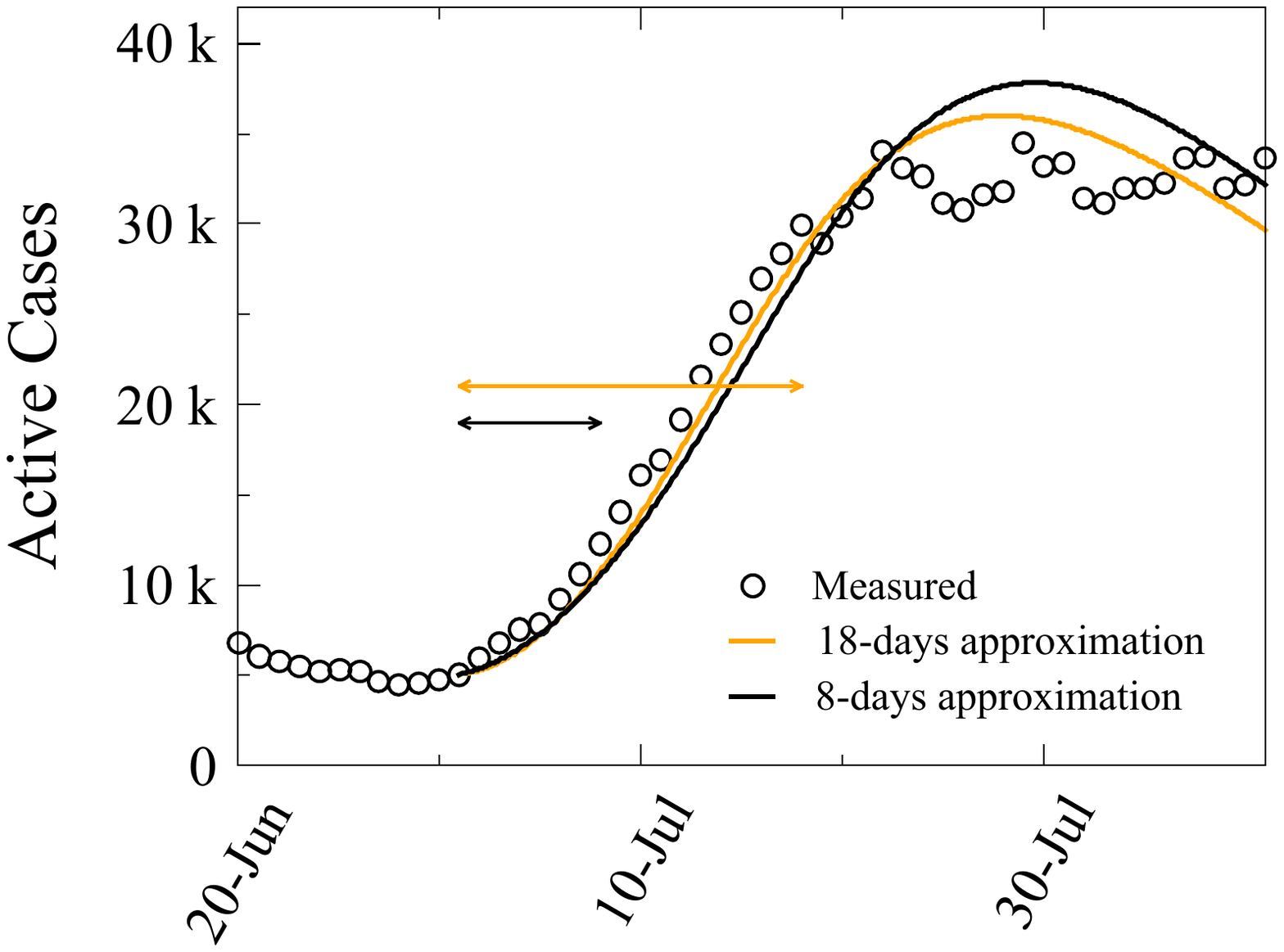}
\vspace{-3 pt}
\caption{The active cases (INSP+ISSP, see Fig.~\ref{fig:1x}) according to the available data (circles) and the solution of the SEInsRD model (curves), fitted to Period A (July 1-8) and Period B (July 1-18) data.~The two periods are shown with horizontal arrows.}
\label{fig:2x}
\end{figure}

The evolution of the active cases from the Worldometer database are shown in Fig.~\ref{fig:2x} with circles (up to 12/8/21).~In addition, Fig.~\ref{fig:2x}  displays the number of active cases, as it was computed from the SEInsRD model, when fitted to Period A (8 days long) and Period B (18 days long) data.

\begin{figure}[t]
\centering
\includegraphics[scale=0.25]{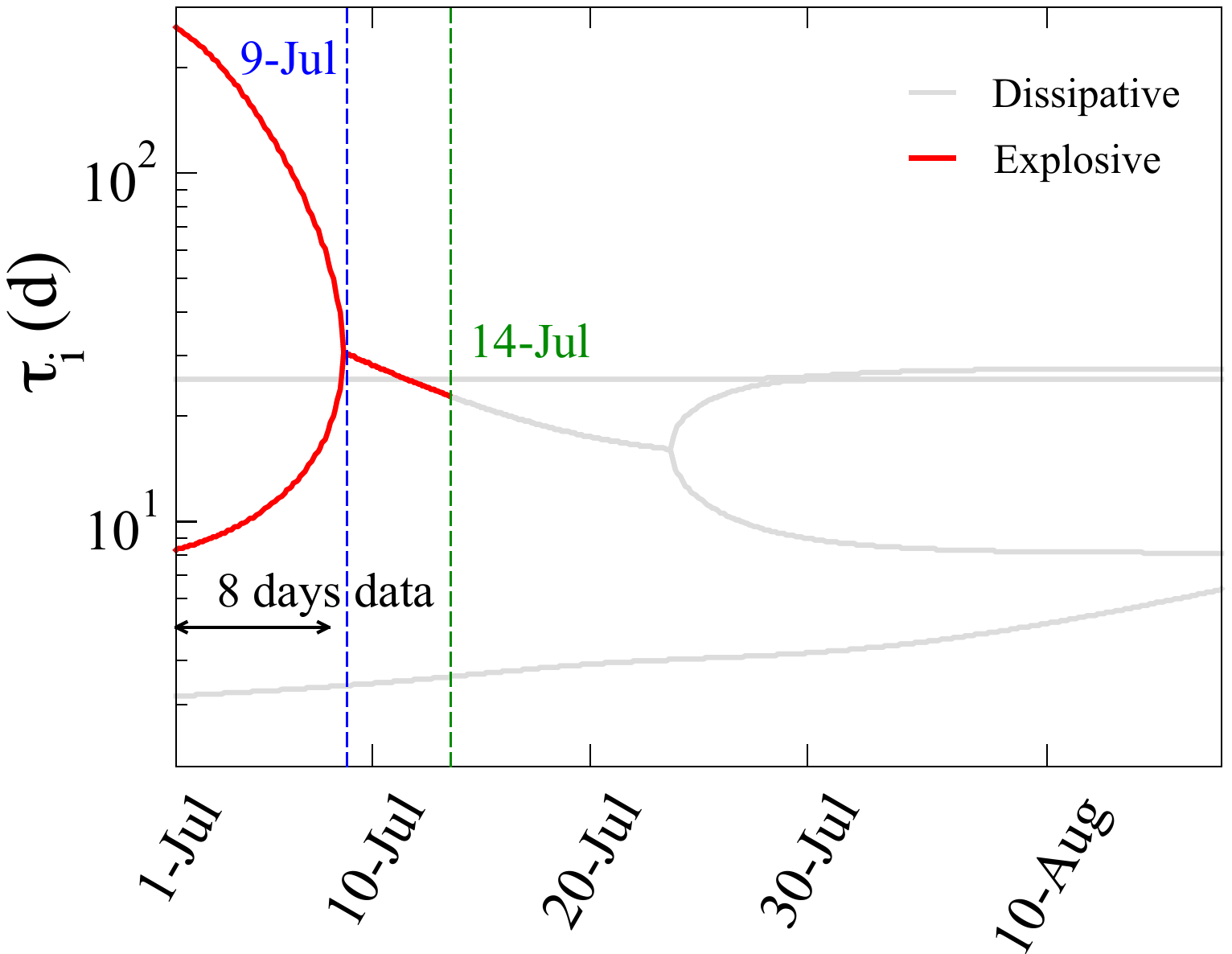}  \hspace{0.15cm} \includegraphics[scale=0.25]{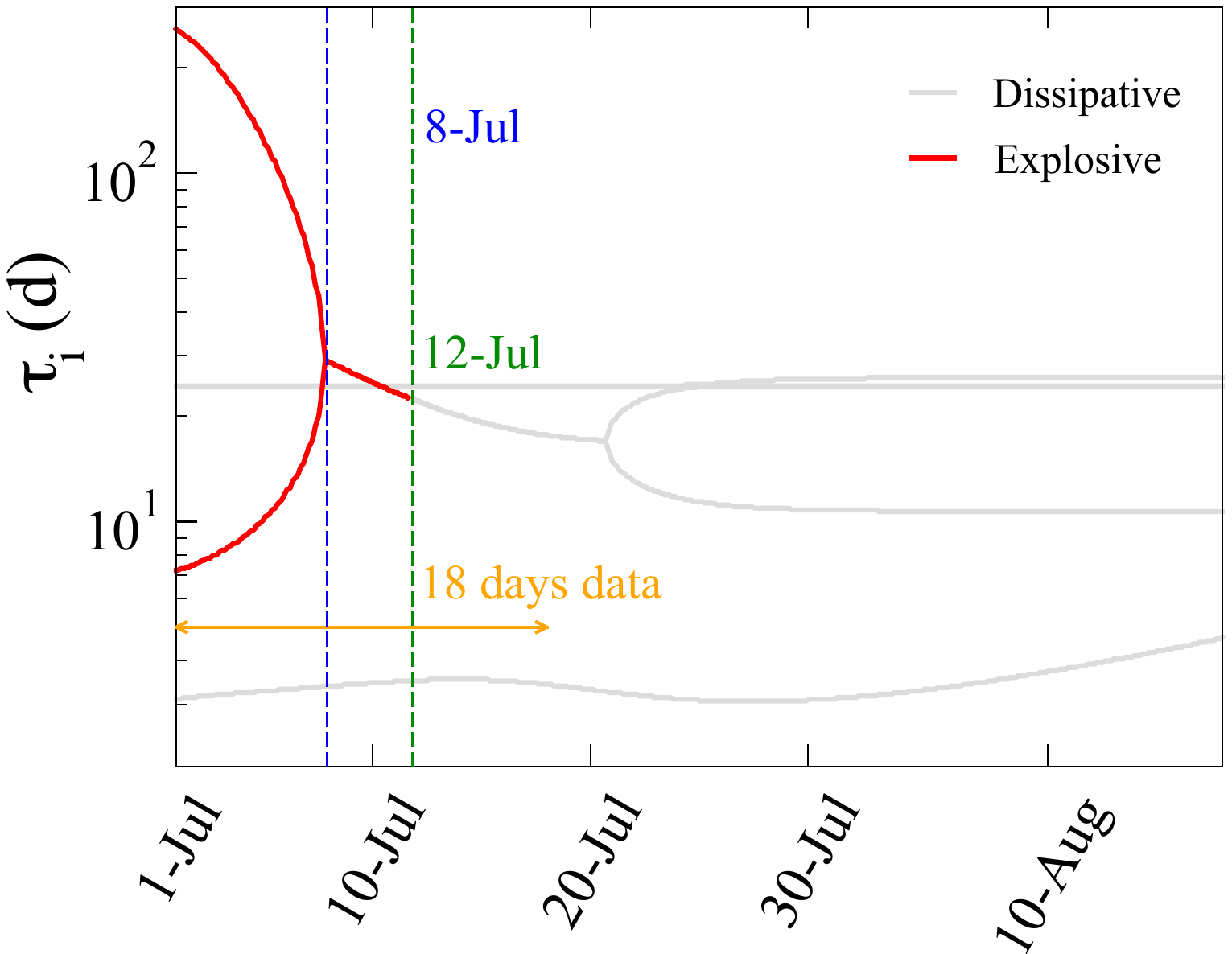}
\vspace{-3 pt}
\caption{The explosive (in red) and dissipative (attenuating, in gray) time scales; profiles based by fitting data from Period A (8 days) and Period B (18 days).~The vertical blue dashed line denotes the point where the two explosive time scales coalesce, indicating that their disappearance is imminent.~The vertical green dashed line denotes the point where the explosive time scales disappear.~The fastest explosive time scale characterises the initial part of the outbreak phase and the fastest dissipative (attenuating) time scale characterises the remaining part.}
\label{fig:3x}
\end{figure}

The  \emph{explosive} and  \emph{dissipative} (attenuating) time scales  that characterise the dynamics in the outbreak phase are shown in Fig.~\ref{fig:3x}.~The \emph{fast explosive} time scale characterises the evolution of the process, from July 1st until the day when this time scale ceases to exist; i.e., July 14 according to Period A and July 12 according to period B.~%When the \emph{explosive} time scales cease to exist, the process is characterised by the fastest implosive time scale.

According to the profiles shown in Fig.~\ref{fig:2x}, the \emph{inflection point} in the profile of the active cases was manifested on July 15 and July 13, according to the solutions based on data spanning Period A and Period B, respectively.~These two days are approximated very well by the days in which the \emph{explosive} time scales disappear; i.e., July 14 and July 12 for the two periods, as shown in Fig.~\ref{fig:3x}.~The disappearance of the \emph{explosive} time scales is preceded by the their coalescence, as shown in Fig.~\ref{fig:3x}; July 9 for Period A and July 8 for Period B.~According to these findings, by July 8-9 it was possible to predict that the development of a peak (plateau in this case) was imminent. 

Given the significance of the \emph{fast explosive} time scale (the faster this scale the more intense the outbreak), it is of interest to identify the paths of the model in Fig.~\ref{fig:1x}, which either promote or oppose its appearance.~The results displayed in Table \ref{tb:TPI_4} suggest that the intensity of the initiation of the outbreak is mainly promoted by Paths 1 and 4 (susceptible becoming exposed and then infected), while it is opposed by Path 6 (infected becoming recovered).~In addition, it is shown that - as the end of the explosive period is reached - the promoting influence of Paths 1 and 4 diminishes, while the opposing one of Path 6 increases (doubles).

\begin{table}[!ht]
\caption{The percentage contribution of the various paths to the development of the fastest explosive time scale, estimated on July 1 and 7. }
\centering
\begin{tabular}{ c || c | c } 
& \multicolumn{2}{c}{Period A solution} \\
\hline
Path & July 1 &  July 7 \\
 \hline \hline
$R^1$		& ~51.8\% &  ~43.8\%  \\
$R^4$		& ~32.6\% &  ~28.9\%  \\
$R^6$		& -14.8\% & -26.3\%  
\end{tabular}
\qquad \qquad
\begin{tabular}{ c || c | c } 
& \multicolumn{2}{c}{Period B solution} \\
\hline
Path & July 1 &  July 7 \\
 \hline \hline
$R^1$		& ~47.6\% &  ~35.1\%  \\
$R^4$		& ~38.3\% &  ~33.4\%  \\
$R^6$		& -13.6\% & -30.3\%  
\end{tabular}
\label{tb:TPI_4}
\end{table}

%\subsection{Conclusions on the 4th wave}
%\label{sec:4th_wave_conclusions}

%It was shown that t
The analysis of available data extending in the period July 1-8 produced similar results with the analysis of data extending in the period July 1-18.~By concentrating on the inflection point in the profile of active cases, it was shown that \textbf{we could know on July 8 or 9 - the latest -} that the peak of active cases was within reach (eventually it was recorded around July 22).~Finally, it was shown that it is possible to \textbf{quantify}, the degree to which the intensity and duration of the \emph{outbreak phase} can be controlled by decreasing the rate by which susceptible become infected and by increasing the rate by which infected become recovered.

%-------------------------------------------------------------------------------------------------------------------------------------------------

\subsection{The 5th wave}
\label{sec:5th_wave}
In this section, results for the 5th wave will be presented.~For validation purposes, the analysis was done by considering solutions obtained by fitting data from two periods in 2021:
\begin{enumerate}[i)]
   \item Period C, October 10 - October 23, 14 days long,
   \item Period D, October 10 - November 6, 28 days long,
\end{enumerate}
where October 10 marks the start of the 5th wave.

\begin{figure}[h!]
\centering
\includegraphics[scale=0.30]{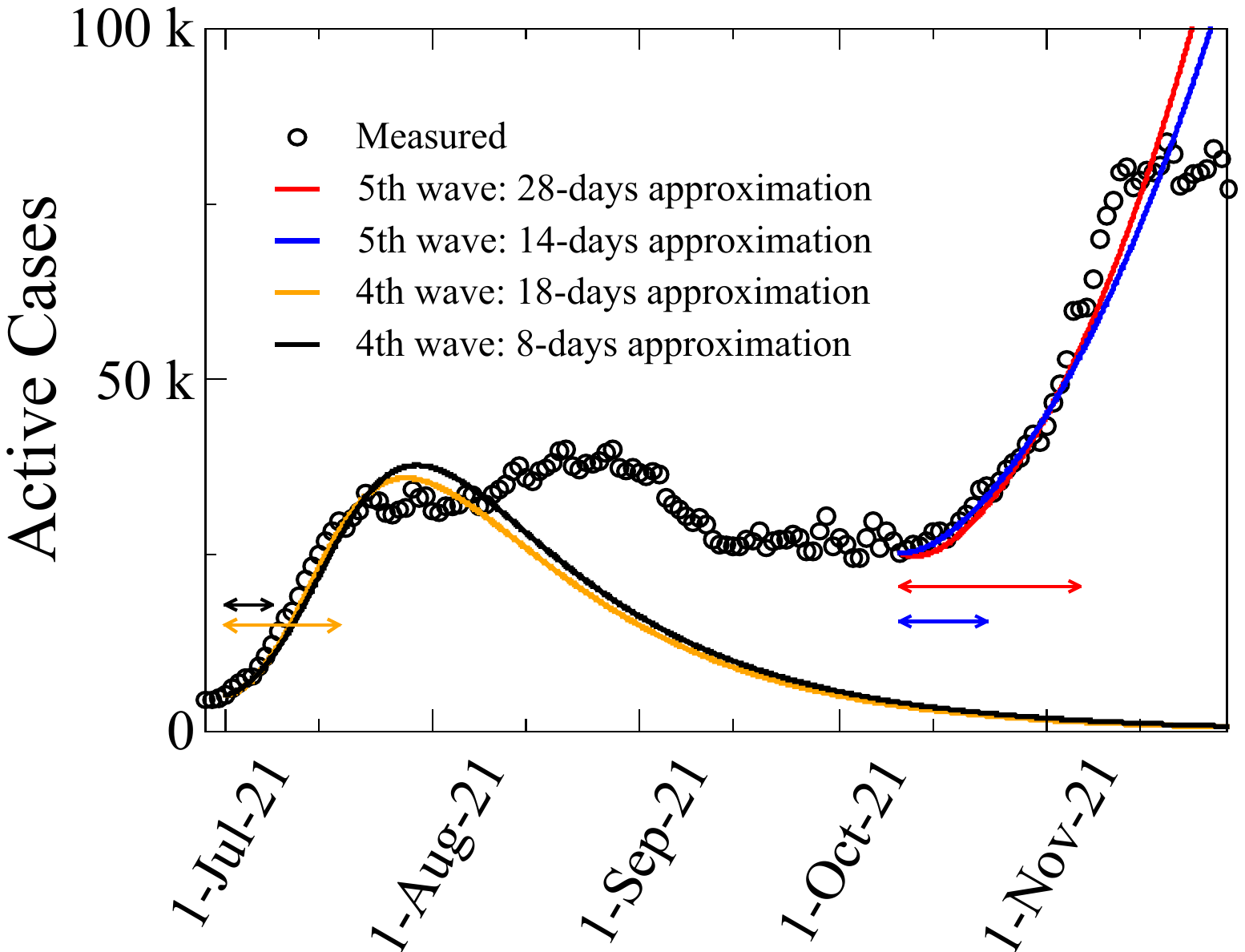}
\vspace{-3 pt}
\caption{The active cases (INSP+ISSP, see Fig.~\ref{fig:1x}) according to the available data (circles) and the solution of the SEInsRD model (curves), fitted to Period C (October 10 - 23) and Period D (October 10 - November 6) data.~The two periods are shown with horizontal arrows.~For comparison, the solution of the model fitted to Periods A and B are displayed.}
\label{fig:4x}
\end{figure}

%\subsection{Results}
The evolution of the active cases from the Worldometer database are shown in Fig.~\ref{fig:4x} with circles (up to 23/11/21).~In addition, Fig.~\ref{fig:4x}  displays the number of active cases, as it was computed from the SEInsRD model, when fitted to Period C (14 days long) and Period D (28 days long) data; the extend of the two periods shown in Fig.~\ref{fig:4x}.

The  \emph{explosive} and  \emph{dissipative} (attenuating) time scales  that characterise the dynamics in the outbreak phase are shown in Fig.~\ref{fig:5x}.~The \emph{fast explosive} time scale is present from the assumed start of the 5th wave in October 10, until 2 January 2022 according to the Period C data and 26 December 2021 according to Period D data.~Since the influence of the fast explosive time scale weakens after it meets the slow one (21-22 November), it is expected that the inflection point of the infected population will manifest early in December.%When the \emph{explosive} time scales cease to exist, the process is characterised by the fastest implosive time scale.

\begin{figure}[h!]
\centering
\includegraphics[scale=0.25]{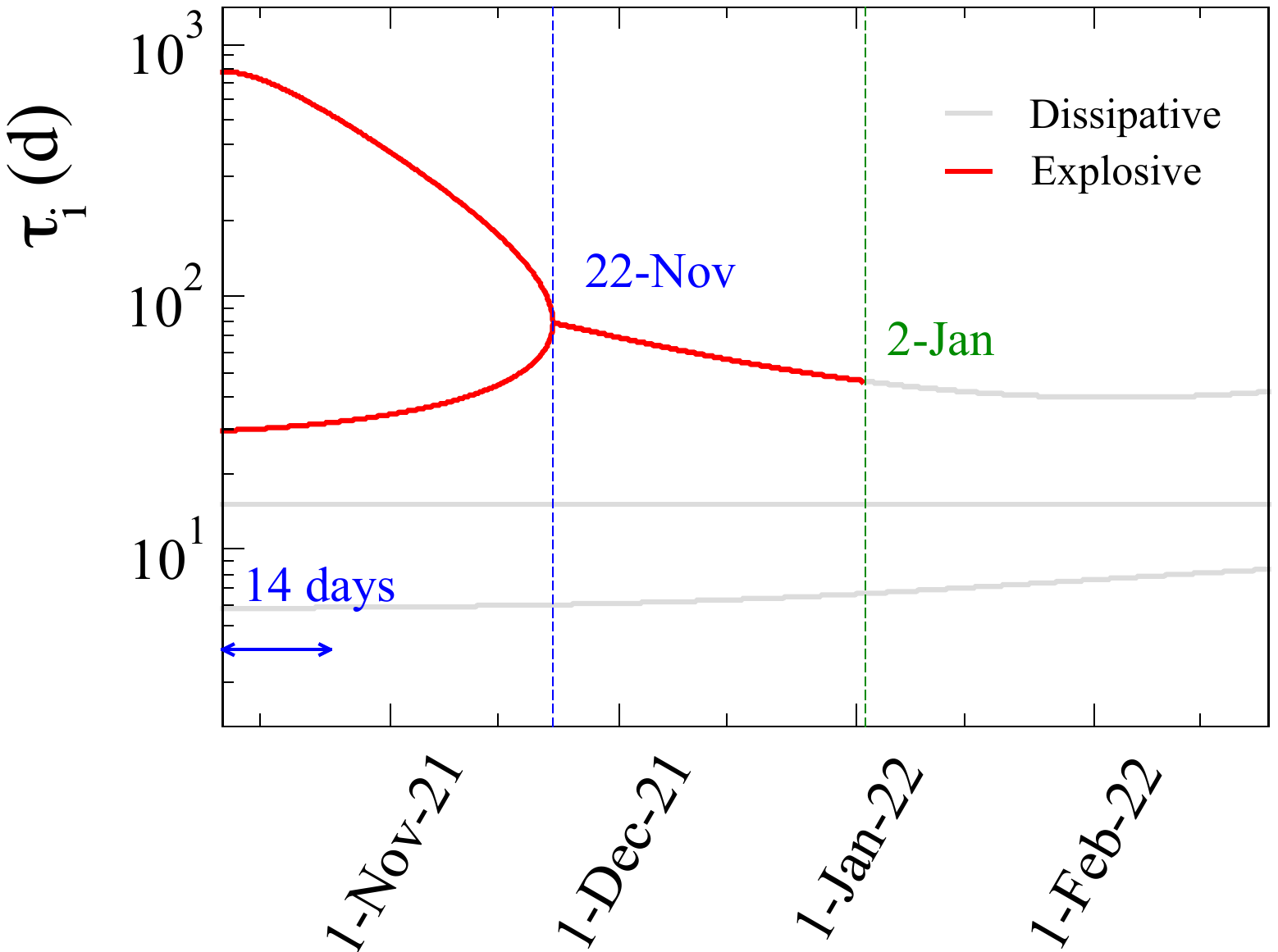}  \hspace{0.15cm} \includegraphics[scale=0.25]{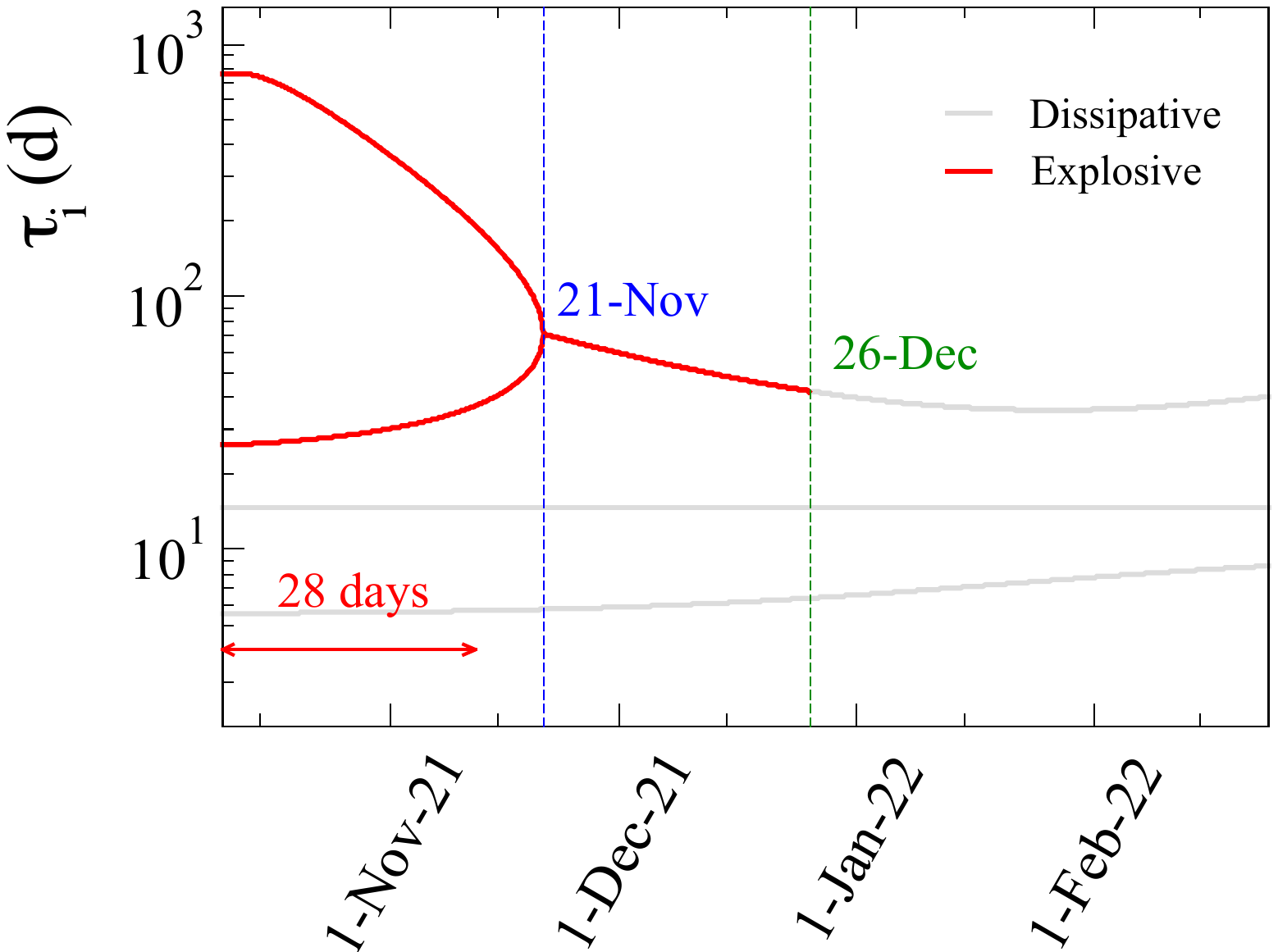}
\vspace{-3 pt}
\caption{The explosive (in red) and dissipative (attenuating, in gray) time scales; profiles based by fitting data from Period C (14 days) and Period D (28 days).~The vertical blue dashed line denotes the point where the two explosive time scales coalesce, indicating that their disappearance is imminent.~The vertical green dashed line denotes the point where the explosive time scales disappear.~The fastest explosive time scale characterises the initial part of the outbreak phase and the fastest dissipative (attenuating) time scale characterises the remaining part.}
\label{fig:5x}
\end{figure}

%-------------------------------------------------------------------------------------------------------------------------------------------------

%\label{sec:5th_wave_results}
Similarly to 4th wave, it is of interest to identify the paths of the model in Fig.~\ref{fig:1x}, which either promote or oppose the generation of the \emph{fast explosive} time scale.~The results displayed in Table \ref{tb:TPI_5} suggest that the intensity of the initiation of this outbreak is mainly promoted by Paths 1 and 4 (susceptible becoming exposed and then infected), while it is opposed by Path 6 (infected becoming recovered).~It is shown that over a 2-weeks period from the start of the 5th wave, there are no significant changes in the major contributions from paths 1, 4 and 6.~This is in contrast to the 4th wave (see Table  \ref{tb:TPI_4}), where over a 1-week period from the start there was a significant decrease in the promoting to the explosive stage contribution of Paths 1 and 4 and a significant increase in the opposing contribution of Path 6.~This feature is related to the expectation that the duration of the 5th wave will be longer than that of the 4th one; see Section \ref{66666}.

In addition, a comparison of Tables \ref{tb:TPI_4} and   \ref{tb:TPI_5} reveals that the opposing to the explosive time scale influence of Path 6 (infected becoming recovered) is much stronger in wave 5 (about -33\% during the first two weeks of the outbreak) than in wave 4 (about -15\% at the start of the wave to about -30\% a week later).

%\textcolor{blue}{The same influence of these paths holds for more that two weeks. In contrast to wave 4, during wave 5, the path where the infected population becomes recovered seems to affect more the dynamics of the system than the path where the exposed population becomes infected.~Even though the model used here does not include variables for the quarantined or vaccined population, the driving mechanism of each wave can be identified by the analysis of the explosive time scale. ---------------------------}

\begin{table}[!ht]
\caption{The percentage contribution of the various paths to the development of the fastest explosive time scale, estimated on October 10 and 24. }
\centering
\begin{tabular}{ c || c | c } 
& \multicolumn{2}{c}{Period C solution} \\
\hline
Path & October 10 &  October 24 \\
 \hline \hline
$R^1$		& ~48.9\% &  ~48.4\%  \\
$R^6$		& -33.4\% &  -34.5\%  \\
$R^4$		& ~17.1\% & ~16.5\%  
\end{tabular}
\qquad \qquad
\begin{tabular}{ c || c | c } 
& \multicolumn{2}{c}{Period D solution} \\
\hline
Path & October 10 &  October 24 \\
 \hline \hline
$R^1$		& ~48.9\% &  ~48.5\%  \\
$R^6$		& -32.2\% &  -33.2\%  \\
$R^4$		& ~18.3\% & ~17.7\%  
\end{tabular}
\label{tb:TPI_5}
\end{table}

%-------------------------------------------------------------------------------------------------------------------------------------------------

%\subsection{Conclusions on the 5th wave}
%\label{sec:5th_wave_conclusions}

%It was shown that 
The analysis of available data extending in the period October 10-23 produced qualitatively similar results with the analysis of data extending in the period October 10 - November 6.~Both periods predict the development of an inflection point (the point where the increase of active cases stops accelerating and starts decelerating) early in December.

\emph{Although the explosive period during the 5th wave is longer than that during the 4th wave, the opposing Path 6 (infected population becoming recovered one) is stronger in the 5th wave from its start, most likely due to the larger portion of vaccinated people in the infected population.~Consequently, if this opposing factor had not being so strong, the promoting influence of Paths 1 and 4 (susceptible population becoming exposed and then infected one) would have made the \emph{outbreak phase} even stronger.~As in the 4th wave, this opposing influence of Path 6 is expected to grow further towards the end of the explosive period.}

%\bibliographystyle{tfq}
%\bibliography{Report_References}

%%%%%%%%%%%%%%%%%%%%%%%%%%%%%%%%%%%%%%%%%%%%%%%%%%%%%%%%%%%%%%%%%%%%%%%%%
\subsection{The 6th wave}
In this section, results for the 6th wave will be presented.~The analysis was done by considering a solution obtained by fitting data from a period of ten days, say Period E, from 26 December 2021 to 4 January 2022, where 26 December 2021 marks the start of the 6th wave.

%\subsection{Results}

\begin{figure}[h!]
\centering
\includegraphics[scale=0.30]{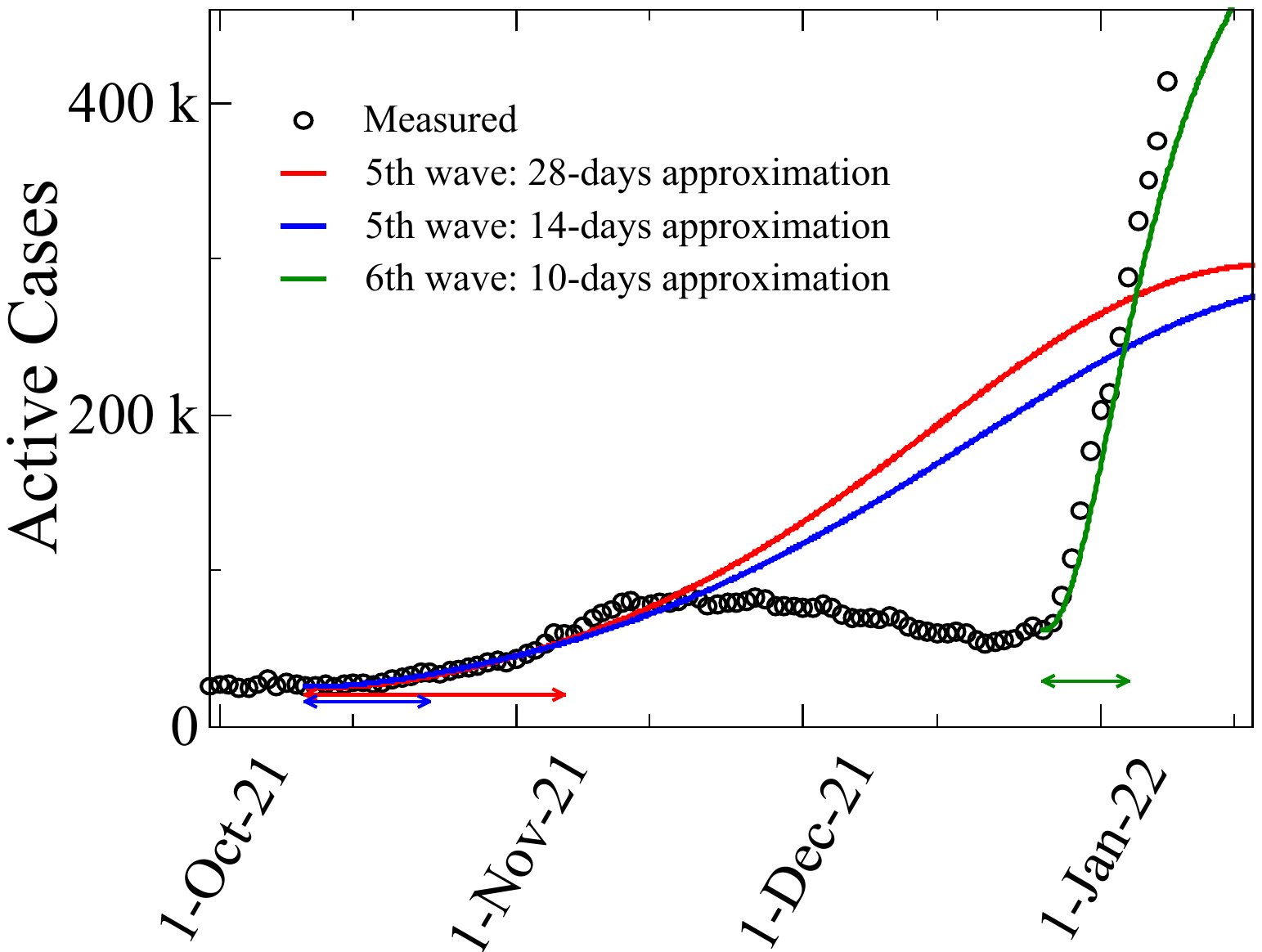}
\vspace{-3 pt}
\caption{The active cases (INSP+ISSP, see Fig.~\ref{fig:1x}) according to the available data (circles) and the solution of the SEInsRD model (curves), fitted to Period E (December 26 '21- January 4 '22) data.~For comparison, the solution of the model fitted to Periods C and D are displayed.~The various periods are shown with horizontal arrows.}
\label{fig:5bx}
\end{figure}

The evolution of the active cases from the Worldometer database are shown in Fig.~\ref{fig:5bx} with circles (up to 08/01/22), along with the number of active cases, as it was computed from the SEInsRD model, when fitted to the Period E data.

The  \emph{explosive} and  \emph{dissipative} (attenuating) time scales  that characterise the dynamics in the outbreak phase are shown in Fig.~\ref{fig:5x}.~The \emph{fast explosive} time scale is present from the assumed start of the 6th wave in December 26 until 31 December 2021.~Since the influence of the fast explosive time scale weakens after it meets the slow one (30 December), it is expected that the inflection point of the infected population will manifest the first days of January 2022.%When the \emph{explosive} time scales cease to 

\begin{figure}[t]
\centering
\includegraphics[scale=0.3]{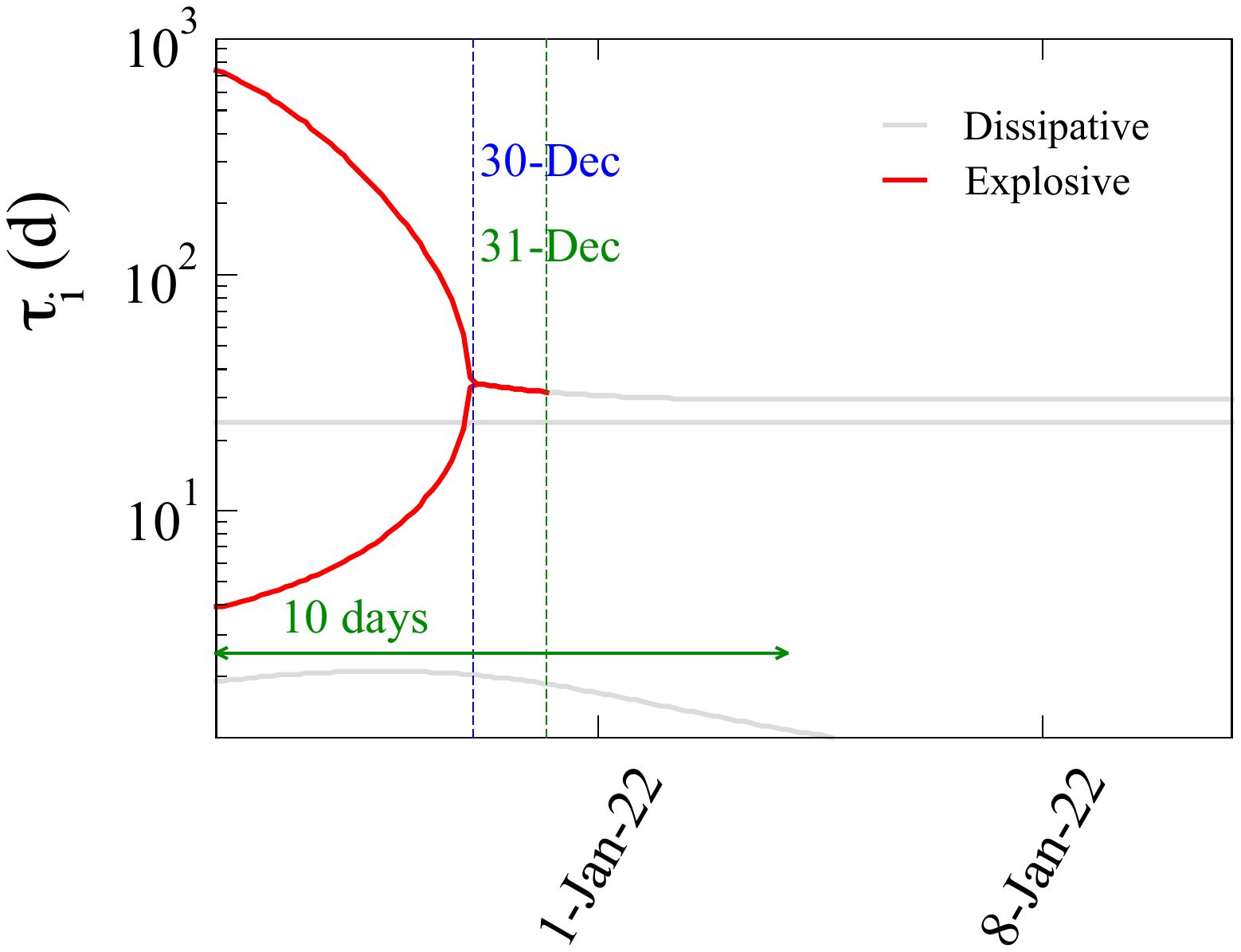}
\vspace{-3 pt}
\caption{The explosive (in red) and dissipative (attenuating, in gray) time scales; profiles based by fitting data from Period E (10 days).~The vertical blue dashed line denotes the point where the two explosive time scales coalesce, indicating that their disappearance is imminent.~The vertical green dashed line denotes the point where the explosive time scales disappear.~The fastest explosive time scale characterises the initial part of the outbreak phase and the fastest implosive (dissipative) time scale characterises the remaining part.}
\label{fig:5x}
\end{figure}

The use of the CSP tools reveals the paths that affect the generation of the \emph{fast explosive} time scale.~The results displayed in Table \ref{tb:TPI_5b} suggest that the intensity of the initiation of this outbreak is mainly promoted by Path 4 (exposed becoming effected) and in a lesser degree by Path 1 (susceptible becoming exposed), while it is opposed by Path 6 (infected becoming recovered).

\begin{table}[!ht]
\caption{The percentage contribution of the various paths to the development of the fastest explosive time scale, estimated on December 26 and 28. }
\centering
\begin{tabular}{ c || c | c } 
& \multicolumn{2}{c}{Period E solution} \\
\hline
Path & December 26 &  December 28 \\
 \hline \hline
$R^4$		& ~56.8\% &  ~59.1\%  \\
$R^1$		& ~35.6\% &  ~28.6\%  \\
$R^6$		& -7.4\% & -11.9\%  
\end{tabular}
\label{tb:TPI_5b}
\end{table}
According the Table~\ref{tb:TPI_5b}, during the outbreak of the 6th wave, the relevant contributions of the main paths do not change much. Furthermore, it is shown that the opposing to the explosive time scale influence of Path 6 (infected becoming recovered) is much more weaker in the 6th wave, compared to the one of both waves 4 and 5; see Tables \ref{tb:TPI_4} and \ref{tb:TPI_5}.

%%%%%%%%%%%%%%%%%%%%%%%%%%%%%%%%%%%%%%%%%%%%%%%%%%%%%%%%%%%%%%%%%%%%%%%%%
%\subsection{Conclusions on the 6th wave}

%It was shown that 
The analysis of available data extending in the period 26 December 2021 to 4 January 2022 predicted the development of an inflection point of active cases (the point where the increase of cases stops accelerating and starts decelerating) early in January.~This implies that the relaxation period (decrease of active cases) will initiate soon.

In contrast to the case of the 4th and 5th waves, the promoting the outbreak phase influence of Path 4 (exposed population becoming infected) was more significant than that of Path 1 (susceptible population becoming exposed), confirming the wider dispersion of the omicron variant to the general population 

It was also shown that the opposing influence of Path 6 (infected population becoming recovered one) is very small; smaller than in the 4th and 5th waves, most likely due to the weaker effects of the vaccines when dealing with the omicron variant.

%%%%%%%%%%%%%%%%%%%%%%%%%%%%%%%%%%%%%%%%%%%%%%%%%%%%%%%%%%%%%%%%%%%%%%%%%
\section{Comparison between waves 4, 5 and 6}
\label{66666}

Tables \ref{tb:6a} and \ref{tb:6b} allow for a direct comparison of the dynamics of the three outbreak phases and the origin of their differences.

Table \ref{tb:6a} lists the explosive time scale $\tau$ at the start of the outbreak phase and at the point where this time scale disappears; the latter coinciding with the inflection point of the active cases.~It is shown that $\tau$ is about 8 d, 28 d and 4 d at the start of the outbreak phase in the case of the 4th, 5th and 6th waves, respectively.~This implies that the process is twice as fast at the start of the outbreak phase of the 6th wave in comparison to that of the 4th wave and seven times as fast in comparison to that of the 5th wave.~Regarding the point where the fast time scale disappears, the 4th and 6th waves are equally fast and about twice as fast in comparison to the 5th wave.

\begin{table}[!ht]
\caption{Comparison between the developing explosive time scales of the 4th, 5th and 6th waves. $\tau_{e}^{0}$ and $\tau_{e}^{end}$ represent the value of the fast explosive time scale at the start of each wave and at the time when the explosive modes disappear, respectively.}
\centering
\begin{tabular}{ cc | cc | c } 

\multicolumn{2}{c}{Wave 4} & \multicolumn{2}{c}{Wave 5} & Wave 6 \\
 \hline  
 Period A & Period B & Period C & Period D & Period E \\
\hline \hline
$\tau_{e,A}^{0}$ = 8.33 d & $\tau_{e,B}^{0}$ = 7.24 d	&$\tau_{e,C}^{0}$ = 29.5 d & $\tau_{e,D}^{0}$ = 25.89 d & $\tau_{e,E}^{0}$ = 3.8 d\\[8pt]
$\tau_{e,A}^{end}$ = 30.53 d & $\tau_{e,B}^{end}$ = 29.14 d	&$\tau_{e,C}^{end}$ = 78.7 d & $\tau_{e,D}^{end}$ = 71.1 d & $\tau_{e,E}^{end}$ = 34.9 d\\
\hline
\end{tabular}
\label{tb:6a}
\end{table}

\begin{table}[!ht]
\caption{Comparison between the major contributors to the explosive time scale at the initiation of the 4th, 5th and 6th waves.}
\centering
\begin{tabular}{ c | c | c } 
Wave 4 & Wave 5 & Wave 6 \\
 \hline  \hline \\[-10pt]
$R^1$:  $\sim$~50\%	& $R^1$:  $\sim$~49\%	&  $R^4$: $\sim$~59\%\\
$R^4$: $\sim$~36\%		& $R^6$:  $\sim$-33\%	&  $R^1$: $\sim$~29\%\\
$R^6$: $\sim$-14\%		& $R^4$:  $\sim$~18\%	& $R^6$: $\sim$-12\%  \\
\hline
\end{tabular}
\label{tb:6b}
\end{table}

Table  \ref{tb:6b} lists the contributions of the various paths of the model employed to the explosive time scale at the start of the outbreak in each wave considered.~The two major findings are (i) the switch from the non-infected population to the exposed one, as the population contributing the most to the outbreak phase, which confirms the wider dispersion of the omicron variant to the general population and (ii) the weaker opposition to the outbreak phase in the 6th wave relative to the 5th one (the first dominated by the omicron variant and the latter by the delta one), which is in accordance to the higher infection rate due to the omicron variant and the decreased efficiency of the vaccines when dealing with this variant..

%================================

%~The initial value of the fast explosive time scale is calculated 3.8 days, which gives roughly an approximation of the duration of the period where the explosive time scales develop and is very short compared to the similar period of the previous waves; i.e. six days.~This indicates that the 6th wave will last shorter and is much steeper in comparison to waves 4 and 5.
  
%~This indicates that the timeframe of action of the explosive time scale that characterises this period where the exposed population becomes effected is much shorter than the timeframe of action of the infected-becoming-recovered process, i.e., the rate which the exposed population becomes effected is much more faster than the one which the infected becomes recovered.

%~Probably, this is due to the larger portion of young or vaccinated people in the infected population.

%%%%%%%%%%%%%%%%%%%%%%%%%%%%%%%%%%%%%%%%%%%%%%%
%%%%%%%%%%%%%%%%%%%%%%%%%%%%%%%%%%%%%%%%%%%%%%%
%%%%%%%%%%%%%%%%%%%%%%%%%%%%%%%%%%%%%%%%%%%%%%%

%%%%%%%%%%%%%%%%%%%%%%%%%%%%%%%%%%%%%%%%%%%%%%%
%%%%%%%%%%%%%%%%%%%%%%%%%%%%%%%%%%%%%%%%%%%%%%%
%%%%%%%%%%%%%%%%%%%%%%%%%%%%%%%%%%%%%%%%%%%%%%%

\section{Conclusions}
\label{Conclusions}

Results based on an alternative and robust methodology for predicting the evolution of each COVID-19 wave were presented.~The methodology focuses on the \emph{fast explosive} time scale that characterises the intensity and duration of the \emph{outbreak phase}.~The point were this time scale ceases to characterize the process coincides with the inflection point of active cases; the point where the increase of active cases stops accelerating and starts decelerating.~Since the inflection point precedes the peak, this methodology serves as an early warning of the peak.~In addition, this methodology allows for the identification of the factors (paths) that can influence the intensity and length of the \emph{outbreak phase}.

The last major three epidemic waves were analyzed, as they were recorded in the periods 1.7.2021-9.10.2021 (4th wave), 10.10.2021-19.12.2021 (5th wave) and 20.12.2021-today (6th wave, still developing); the first two periods dominated by the \emph{delta} variant of the virus and the last one by the \emph{omicron} variant.~By examining the dynamics of their \emph{outbreak phase}, it is concluded that:
\begin{enumerate}[i)]
    \item The most intense \emph{outbreak phase} is that of the 6th wave, followed by that of the 4th wave; the one of the 5th wave being by far the least intense.~This finding confirms (i)  \emph{the dominance of the omicron variant, relative to the delta one}, (ii) \emph{the relaxing influence of vaccinations} (much higher percentage of vaccinated population during the 5th wave relative to the 4th and not significant difference between the 5th and the 6th waves) and (iii) the weak influence of vaccines when dealing with the omicron variant.
    \item The major contribution to the \emph{outbreak phase} of the 4th and 5th waves originates from the non-infected population becoming exposed (infected but not yet able to transmit the virus), followed by the exposed population becoming infected (able to infect).~In the 6th wave, this ordering is reversed, so that the population group dominating the \emph{outbreak phase} is the exposed population and not the non-infected ones of the 4th and 5th waves.~This switch \emph{confirms the wider dispersion of the omicron variant in the general population}.
    \item The major opposing influence to the intensity of the \emph{outbreak phase} of all three waves originates from the infected population becoming recovered.~This opposition is  strongest in the 5th wave (least intense wave) and weakest in the 6th wave (most intense wave).
    \begin{enumerate}[a)]
    \item The stronger opposition to the \emph{outbreak phase}  in the 5th wave relative to the 4th one (both dominated by the delta variant) is \emph{in accordance to the higher vaccinated percentage of the general population during the 5th wave}.
    \item The weaker opposition to the \emph{outbreak phase}  in the 6th wave relative to the 5th one (the first dominated by the omicron variant and the latter by the delta one) is \emph{in accordance to the higher infection rate due to the omicron variant and the decreased efficiency of the vaccines when dealing with this variant}.
    \end{enumerate}
    \item The inflection point of active cases in the current 6th wave occurred in the first days of January 2022, so the outbreak is about to start relaxing.
\end{enumerate}

\textbf{NOTE 1:} \emph{The predictions reported here for the 4th wave were based on the prevailing conditions in Periods A (July 1-8) and B (July 1-18).~It was demonstrated that it was possible to predict by July 8 or 9 - the latest - that the peak of active cases was within reach; eventually, it was recorded around July 22.}

\emph{The predictions reported for the 5th wave are valid under the assumption that the prevailing conditions in Periods C (October 10-23) and D (October 10 - November 6) keep prevailing in the period that follows.}

\emph{Similarly, the predictions reported for the 6th wave are valid under the assumption that the prevailing conditions in Periods E (26 December 2021 to 4 January 2022) keep prevailing in the period that follows.}

~~

\textbf{NOTE 2:} \emph{In accordance to recent studies, it is assumed here that the probability of transmission between vaccinated or not-vaccinated does not differ significantly}.

%%%%%%%%%%%%%%%%%%%%%%%%%%%%%%%%%%%%%%%%%%%%%%%%%%
%%%%%%%%%%%%%%%%%%%%%%%%%%%%%%%%%%%%%%%%%%%%%%%%%%
%\section*{Acknowledgments}
%\label{Acknowledgments}
%This work was sponsored by competitive research funding from King Abdullah University of Science and Technology (KAUST). The support from Khalifa University of Science and Technology, via project 8434000269, is gratefully acknowledged.
%%%%%%%%%%%%%%%%%%%%%%%%%%%%%%%%%%%%%%%%%%%%%%%%%%
%%%%%%%%%%%%%%%%%%%%%%%%%%%%%%%%%%%%%%%%%%%%%%%%%%
%This work was sponsored by competitive research funding from King Abdullah University of Science and Technology (KAUST).%~The simulations presented in this work made use of KAUST Supercomputing Laboratory

%% References can be added with or without bibTeX database
%%
%% References with bibTeX database:
%% The bibliography style, elsarticle-num.bst, is used and available within the template package

\appendix
\section{Appendix}\label{sec:appendix}
According to the frequency dependent representation of SEInsRD model in Eq.~\eqref{eq:SEI2RD}, the following conclusions can be drawn, regarding to the inflection points $SP$, $EP$, $INSP$ and $ISSP$, through analytic calculations.

%%%%%%%%%%%%%%%%%%%%%%
\subsection{Inflection point of SP close to max R$^1$ and R$^2$}

Differentiation of the SP differential equation in Eq.~\eqref{eq:SEI2RD} yields: 
\begin{equation}
\ddot{SP}=-\dot{R}^1-\dot{R}^2-\dot{R}^3=-(\dot{R}^1+\dot{R}^2)+\mu_{TP} (R^1+R^2+R^3)
\end{equation}

Neglecting the terms $\mu_{TP} (R^1+R^2+R^3)$, since $\mu_{TP}\ll1$, yields:
\begin{equation}
\ddot{SP}\approx-(\dot{R}^1+\dot{R}^2)
\label{eq:SPdd}
\end{equation}

During the explosive stage both $R^1$ and R$^2$ increase with time; i.e. initially the rate by which $SP$ gets infected increases with time ($\dot{SP}<0$ and $\ddot{SP}<0$).~According to Eq.~\eqref{eq:SPdd}, an inflection point of $SP$ ($
\ddot{SP}=0$) might occur when R$^1$ and R$^2$ attain their maximum values ($\dot{R}^1=\dot{R}^2=0$). The inflection point of SP and the points in time where R$^1$ and R$^2$ attain their maximum values are indeed very close.~The differences are related to the magnitude of the term $\mu_{TP} (R^1+R^2+R^3)$.

%%%%%%%%%%%%%%%%%%%%%%
\subsection{Inflection points of SP and EP}

Differentiation of the sum of SP and EP equations in Eq.~\eqref{eq:SEI2RD} yields: 
\begin{equation}
\ddot{SP}+\ddot{EP}=-\dot{R}^3-\dot{R}^4-\dot{R}^5=-\dfrac{\dot{EP}}{aincp}-\mu_{TP} (\dot{SP}+\dot{EP})
\end{equation}

Neglecting the terms $\mu_{TP} (\dot{SP}+\dot{EP})$, since $\mu_{TP}\ll1$, yields: 
\begin{equation}
\ddot{SP}+\ddot{EP}\approx-\dfrac{1}{aincp} \dot{EP}
\label{eq:SP+EPdd}
\end{equation}

Given that $\dot{EP}>0$ during the explosive stage, Eq.~\eqref{eq:SP+EPdd} indicates that $\ddot{SP}+\ddot{EP} < 0$ there.~In addition, with the exception of the very early and very late parts of the explosive stage, $\ddot{SP}\le 0$ and $\ddot{EP} \ge 0$ throughout this stage; i.e., the rate of increase of $EP$ and the rate of decrease of $SP$ intensify.~Moving towards the end of the explosive stage, $\dot{EP}$ decreases, being reflected with a drop of the increasing trend of EP.~At the end of the explosive stage, $\dot{EP}$ attains small positive values, so that Eq.~\eqref{eq:SP+EPdd} yields that the inflection point of EP ($\ddot{EP}=0$) comes shortly before the end of the explosive stage, followed by the inflection point of SP ($\ddot{SP}=0$).

%%%%%%%%%%%%%%%%%%%%%%
\subsection{Inflection points of INSP and ISSP coinciding with max R$^4$}

Differentiation of the sum of INSP and ISSP equations in Eq.~\eqref{eq:SEI2RD} yields:
\begin{align}
\ddot{INSP}+\ddot{ISSP}& =\dot{R}^4-\dot{R}^6-\dot{R}^7-\dot{R}^8-\dot{R}^9 \nonumber \\
& =\dfrac{\dot{EP}}{aincp}-\left( \dfrac{1}{aip} +\mu_{INSP} \right) \dot{INSP}-\left( \dfrac{1}{aip} +\mu_{ISSP} \right) \dot{ISSP}
\label{eq:IN+ISInfl}
\end{align}

The summation of all the differential equations in Eq.~\eqref{eq:SEI2RD} yields:
\begin{equation}
\dot{SP}+\dot{EP}+\dot{INSP}+\dot{ISSP}+\dot{RP}+\dot{IDP}=-R^3-R^5-R^{10} = 
\end{equation}
in which the terms $\dot{RP}$, $\dot{IDP}$ and $-R^3-R^5-R^{10}$ are very small during the explosive stage, such that:
\begin{equation}
\dot{SP}+\dot{EP}+\dot{INSP}+\dot{ISSP}=0
\label{eq:temp}
\end{equation}
Differentiation of Eq.~\eqref{eq:temp} yields:
\begin{equation}
\ddot{SP}+\ddot{EP}+\ddot{INSP}+\ddot{ISSP}=0
\end{equation}
which, after algebraic calculations, takes the form:
\begin{equation}
\mu_{TP} (\dot{SP}+\dot{EP})+\left( \dfrac{1}{aip}+\mu_{INSP} \right)\dot{INSP}+\left( \dfrac{1}{aip}+\mu_{ISSP} \right)\dot{ISSP}=0
\label{eq:temp1}
\end{equation}
Differentiation of Eq.~\eqref{eq:temp1} implies:
\begin{equation}
\mu_{TP} (\ddot{SP}+\ddot{EP})+\left( \dfrac{1}{aip}+\mu_{INSP} \right)\ddot{INSP}+\left( \dfrac{1}{aip}+\mu_{ISSP} \right)\ddot{ISSP}=0
\end{equation}
in which the substitution of Eq.~\eqref{eq:SP+EPdd} yields:
\begin{equation}
-\left( \mu_{TP}^2+\dfrac{\mu_{TP}}{aincp}  \right)\dot{EP}+\left( \dfrac{1}{aip}+\mu_{INSP} \right)\ddot{INSP}+\left( \dfrac{1}{aip}+\mu_{ISSP} \right)\ddot{ISSP}=0
\label{eq:IN+ISInfl}
\end{equation}

Given that $\dot{EP}>0$, $\ddot{INSP}>0$ and $\ddot{ISSP}>0$ during the explosive stage and that $\dot{INSP}$ and $\dot{ISSP}$ attain small positive values, when moving towards to the end of it, Eq.~\eqref{eq:IN+ISInfl} denotes that the inflection points of INSP and ISSP ($\ddot{INSP} \approx \ddot{ISSP}=0$) coincide with the time point at which EP attains its maximum value ($\dot{EP}=0$).
\bibliography{covid_biblio}% common bib file
\bibliographystyle{unsrtnat}

%%% Uncomment this section and comment out the \bibliography{references} line above to use inline references.
% \begin{thebibliography}{1}

% 	\bibitem{kour2014real}
% 	George Kour and Raid Saabne.
% 	\newblock Real-time segmentation of on-line handwritten arabic script.
% 	\newblock In {\em Frontiers in Handwriting Recognition (ICFHR), 2014 14th
% 			International Conference on}, pages 417--422. IEEE, 2014.

% 	\bibitem{kour2014fast}
% 	George Kour and Raid Saabne.
% 	\newblock Fast classification of handwritten on-line arabic characters.
% 	\newblock In {\em Soft Computing and Pattern Recognition (SoCPaR), 2014 6th
% 			International Conference of}, pages 312--318. IEEE, 2014.

% 	\bibitem{hadash2018estimate}
% 	Guy Hadash, Einat Kermany, Boaz Carmeli, Ofer Lavi, George Kour, and Alon
% 	Jacovi.
% 	\newblock Estimate and replace: A novel approach to integrating deep neural
% 	networks with existing applications.
% 	\newblock {\em arXiv preprint arXiv:1804.09028}, 2018.

% \end{thebibliography}

\end{document}